# ON SOME CONDITIONS ON A NOETHERIAN RING.

C.L.Wangneo

Jammu,J&K,India,180002

(E-mail:-wangneo.chaman@gmail.com )

**Abstract :-** In this paper for a noetherian ring R with nilradical N(R) we define semiprime ideals P(R) and Q(R) called as the left and right krull homogenous parts of N(R) .We also recall the known definitions of localisability and the weak ideal invariance (w.i.i for short ) of an ideal of a noetherian ring R . We then state and prove results that culminate in our main theorem whose statement is given below ;

**<u>Theorem</u> :-** Let R be a noetherian ring with nilradical N . Let P(R) and Q(R) be semiprime ideals of R that are the right and left krull homogenous parts of N respectively . Then the following conditions are equivalent ;

(i) N(R) is a right w.i.i ideal of R ( respectively N(R) is a left w.i.i ideal of R ) .

(ii) P(R) is a right localizable ideal of R ( respectively Q(R) is a left localizable ideal of R ) .

## <u>Introduction</u> :-

This paper is divided into two sections . In the brief section (1) we first recall the preliminary definitions and results that are needed to state and prove our main theorem in section (2) .Thus in section (1) we recall from [6] the definitions of ideals A(R) , B(R) of R and the subsets Λ(R) , Λ'(R) of spec.(R) . We also recall the definitions of semiprime ideals P(R) and Q(R) of R called as the right and left krull homogenous parts of the nilradical N(R) respectively . We then state some basic results in this section that are needed to prove our main theorem in section (2) .

Next in section (2) we recall the definitions of the weak ideal invariance and localisability of an ideal of a noetherian ring R as stated for example in [ 1 ]. We then state results in this section that lead us to prove our main theorem namely , theorem (2.10) , that relates right and left w.i.i of the nilradical N(R) of a noetherian ring R to the right and left localisability of the ideals P(R) and Q(R) of R respectively .

**<u>Notation and Terminology</u>:-**

We adhere throughout in this paper to the same terminology as in [6] which is the same as in [2] . Thus a ring is meant to be an associative ring with identity which is not necessarily a commutative ring . By a noetherian ring R we mean that R is both a left as well as a right noetherian ring. By a module M over a ring R we mean that M is a right R-module unless stated otherwise. For the basic definitions regarding noetherian modules and noetherian rings and all those regarding Krull dimension, we refer the reader to [2]. For a noetherian ring B , if M is a right noetherian module over the right noetherian ring B , then we will denote by $|M|_r$ the right krull dimension of M . Similarly $|M|_l$ denotes the left krull dimension of M in case M is a left noetherian B -module over the noetherian ring B . Recall from [2] , a module M is said to be right ( respectively left ) Krull-homogenous if every non zero right (respectively left ) sub -module N of M has same Krull dimension as that of M . Also for the definition of the weak ideal invariance (w.i.i for short) and localizability of an ideal we refer the reader to [1]. We now present some more symbols that we shall use throughout.

**Further Notation** :- If R is a ring then we denote by Spec.R, the set of prime ideals of R and by min.Spec.R, the set of minimal prime ideals of R. Moreover for a right R module M over a ring R we denote by Ass.M the set of prime ideals of R associated to M on the right. A right (or a left) R-module M is said to be right (or left) primary if it has a unique right (or left) associated prime ideal. Also for a ring R we denote by r(T) the right annihilator of a subset T of a right R mdule M and similarly l(T) denotes the left annihilator of a subset T of W in case W is a left R module. For two subsets A and B of a given set, A ≤ B means B contains A and A<B denotes A ≤ B but A≠B. Also A⊄B denotes the subset B does not contain the subset A. For an ideal A of R c(A) denotes the set of elements of R that are regular modulo the ideal A.

**Section (1) (Preliminaries)** :- In this section we first recall the following definitions from [6].

**Definition and Notation(1.1)** :- Let R be a noetherian ring. Denote by A(R) the sum of all right Ideal I of R such that $|I|_r < |R|_r$. Similarly denote by B(R) the sum of all left ideal K of R with $|K|_l < |R|_l$. If there is no confusion regarding the underlying ring R we may write A and B for A(R) and B(R) respectively.

We now state the following basic result regarding the above definition from [6];

**Proposition (1.2)** :- Let R be a Noetherian ring. Let A and B be as in definition(1.1) above. Then the following holds true ;

A is an ideal of R and is the unique largest right ideal of R with $|A|_r < |R|_r$ and $|R/A|_r = |R|_r$. Moreover R/A is a right k-homogenous ring and if $A \neq 0$, then $r(A) \neq 0$ and $|R/r(A)|_r < |R|_r$.

**Proof** :- For the proof of the above result we ask the reader to see the proof of the proposition (1.3) of [6].

We now state proposition(1.3) below which is similar to proposition (1.2) above and may be considered as its left analogue ;

**Proposition(1.3)** :- Let R be a Noetherian ring. Let A and B be as in definition(1.1) above. Then the following holds true ;

B is an ideal of R and is the unique largest left ideal of R with $|B|_l < |R|_l$ and $|R/B|_l = |R|_l$. Moreover R/B is a left k-homogenous ring and if $B \neq 0$, then $l(B) \neq 0$ and $|R/l(B)|_l < |R|_l$.

Next we recall the following definition and notation from [6].

**Definition and Notation (1.4)** :- **(i)** Let R be a Noetherian ring with nilradical R. Define subsets of spec. R denoted by $\Lambda(R)$ and $\Lambda'(R)$ as follows ;

$\Lambda(R) = \{ P_i \text{ in spec. } (R) \ / \ |R/P_i|_r = |R|_r \}$ and

$\Lambda'(R) = \{ Q_i \text{ in spec. } (R) \ / \ |R/Q_i|_l = |R|_l \}$.

If there is no confusion about the underlined ring R we set $\Lambda = \Lambda(R)$ and $\Lambda' = \Lambda'(R)$.

**(ii)** It is not difficult to see that if , $P_i \varepsilon \Lambda(R)$, then $P_i \varepsilon$ min. spec.(R) . We denote by P(R) the ideal P(R) = $\cap P_i$ , $P_i \varepsilon \Lambda$ and by Q(R) the ideal Q(R) = $\cap Q_i$ , $Q_i \varepsilon \Lambda'$ . Then the ideals P(R) and Q(R) are called the right and the left krull homogenous parts respectively of the nilradical N(R) of the ring R . If there is no confusion regarding the underlying ring R we may write P for P(R) and Q for Q(R) . We may note that in this case the factor rings R/P and R/Q are semiprime right and left krull homogenous rings respectively . It is also clear that $N \leq P$ and $N \leq Q$ .

**(iii)** For an ideal I of R we have ; $\Lambda(R/I) = \{ P_i/I$ in spec. (R/I) / $I \leq P_i$ and $|R/P_i|_r = |R/I|_r \}$ and $\Lambda'(R/I) = \{ Q_i/I$ in spec. (R/I) / $I \leq Q_i$ and $|R/Q_i|_l = |R/I|_l \}$ .

**Remark (1)** :- For an ideal I of a noetherian ring R with $|R/I|_r = |R|_r$, it is clear that if N(R) is the nilradical of R , then $N(R/I) = \cap q_i /I$ , where $q_i$ are those prime ideals of R such that $I \leq q_i$ , and qi are the prime ideals minimal over the ideal I ( qi need not be minimal prime ideals of R ) . However $P(R/I) = \cap p_i /I$, where $p_i$ are those prime ideals of R such that $p_i \varepsilon \Lambda$ and $I \leq p_i$ ( pi are necessarily the minimal prime deals of R ) .

**Section (2)** :- Now Following [1] we recall the definition of the weak ideal invariance ( w.i.i for short ) of an ideal of a noetherian ring R . Then , for a noetherian ring R , if P is the semiprime ideal , P =$\cap$p , p $\varepsilon \Lambda$ , we show that the right w.i.i of the nilradical N of R is equivalent to the right w.i.i of the ideal P of R . Note that we have called P as the right krull homogenous part of the nilradical N of the ring R.

**Defiinition (2.1) (Weak ideal invariance)** :- Let R be a noetherian ring. If I is any ideal of R we call I a right weakly ideal invariant ideal ( right w.i.i for short ) if for any right ideal J of R with $|R/J|_r < |R/I|_r$, we have that $|I/JI|_r < |R/I|_r$. Similarly we define the notion of the left w.i.i of an ideal I of a noetherian ring R. An ideal is said to be w.i.i if it is left w.i.i as well as right w.i.i.

We now state for our purposes the following result regarding w.i.i of an ideal as below ;

**Lemma (2.2) :-** Let R be a right noetherian ring with $|R|_r = a$. Then the following hold true ;
**(a)** If R is semiprime then any minimal prime ideal of R is right w.i.i .
**(b)** If an ideal T of R with $|R/T|_r = |R|_r = a$ is right weakly ideal invariant (w.i.i for short), then for a finitely generated module M, if M has an essential Krull- homogenous sub-module H, with $|H|_r = a$ such that HT=0 and $|M/H|_r < |H|_r$, then MT=0.

**Proof** :- **(a)** For the proof of (a) see lemma (2.4) of [1].
**(b)** For the proof of (b) see [5] lemma (2).

**Remark (2)** :- For a semiprime noetherian ring R the following is always true ; For a finitely generated critical module M if Ass.(M) = Q, for some minimal prime ideal Q of R, then $|M|_r = |R/Q|_r$ implies that MQ =0.

**Proof :-** We leave the proof of this fact to the reader.

We state Theorem (2.3) below which gives conditions equivalent to the w.i.i of the nilradical of a noetherian ring R.

**Theorem (2.3)** :- Let R be a noetherian ring with nilradical N. Let A, B and $\Lambda$, $\Lambda'$ be as defined in section (1) above. Let $P = \cap P_i$, $P_i \in \Lambda$. Then the following conditions are equivalent;

**(i)** N is right w.i.i ideal.

**(ii)** For any prime ideal Q of R with $|R/Q|_r = |R|_r$, Q is right w.i.i.

**(iii)** For any ideal I of R with $|R/I|_r = |R|_r$, I is right w.i.i.

**(iv)** If M is any finitely generated critical right R-module with $|M|_r = |R|_r$ and Ass.(M) = Q, then r(M) = Q.

**(v)** P is right w.i.i.

**Proof** :- For the proof of equivalence of (i), (ii), (iii) and (iv) we use theorem (2.3) and theorem (2.5) of [1].

The proof of (iii) => (v) is obvious because $|R/P|_r = |R|_r$.

We now prove (v)=>(iv). To prove it assume (v) and let M be any f.g critical critical right R-module with $|M|_r = |R|_r$ with Ass.(M) = Q. We show if r(M) = I, then I = Q. Let N be the largest submodule of M such that NQ = 0. Now observe that since $|R/Q|_r = |R|_r$, hence we have that $P \leq Q$. It then follows that NP = o. Applying lemma (2.2) (b) we get that MP = 0. Thus M is a R/P- module. Hence from lemma(2.2) (a) Q/P is a right w.i.i ideal of the semiprime ring R/P. Applying lemma (2.2)(b) again implies that we must have MQ = 0. This proves (iv).

Hence the proof of the theorem is complete.

We now recall the definition of localizability of an ideal from [1].

**Definition (2.4)** **( localisability ) :-** Let I be an ideal of a ring R. Let $C(I) = \{c \in R / c+I \text{ is regular in } R/I\}$. The ideal I is said to be right localisable if the elements of $C(I)$ satisfy the right ore condition; that is, given $c \in C(I)$ and $x \in R$, there exists $d \in C(I)$ and $y \in R$ such that $xd = cy$. An ideal is said to be localisable if it is left and right localisable. Define $T(I) = \{x \in R / xc = 0, \text{ for some } c \in C(I)\}$. Note that, if I is right localisable, then $T(I)$ is actually an ideal.

We now address the important question of the localisability of the semiprime ideal $P = \cap P_i$, $P_i \in \Lambda$, of a noetherian ring R in terms of the w.i.i of P.

We first state the following result below that is got by Combining proposition (3.1) and lemma (3.6) or theorem (3.8) of [1];

**Theorem (2.5)** :- Let P be a semiprime ideal of a noetherian ring R such that R/P is right krull homogenous. Then the following hold true;

**(i)** P is right localisable implies that P is right w.i.i.

**(ii)** If P is right w.i.i and P has the right AR property then P is right localisable.

Before stating the key theorem of this section we prove the following results below;

**Theorem (2.6) :-** Let R be a noetherian ring with nilradical N. Let A, B and $\Lambda$ be as defined in section (1) above. Let $P = \cap P_i$, $P_i \in \Lambda$. Then the following are equivalent;

**(a)** N is right w.i.i.

**(b) P** is right w.i.i.

**(c)** For any factor ring R/I of R with $|R/I|_r = |R|_r$, the nilradical N(R/I) of R/I is right w.i.i .

**(d)** N(R/A) is right w.i.i .

**Proof:-** The equivalence of (a) , (b) and (c) follows from Theorem (2.3) above .

(c) => (d) follows by applying theorem (2.3) and by observing from theorem (1.2) that $|R/A|_r = |R|_r$ .

**(d) =>(a)** :- To see this let M be a finitely generated critical R-module with $|M|_r = |R|_r$ and such that Ass.(M) = p , for some prime ideal p of R . We show Mp =0 . Observe that since Ar(A) =0 , and since by proposition (1.2) above , $|R/r(A)| < |R|$ , so we must have that MA=0. Thus M is an R/A - module . Now by the given condition (d) N(R/A) is right w.i.i , hence using theorem (2.3) we must have that Mp=0 . Thus theorem (2.3) guarantees that N is right w.i.i . Before we proceed further we must state the left analogue of theorem (2.6) above .

**Theorem (2.7) :-** Let R be a noetherian ring with nilradical N . Let A , B and Λ' be as defined in section (1) above . Let Q = ∩ $Q_i$ , $Q_i$ ε Λ' . Then the following are equivalent ;

**(a)** N is left w.i.i .

**(b)** Q is left w.i.i **.**

**(c)** For any factor ring R/I of R with $|R/I|_l = |R|_l$, the nilradical N(R/I) of R/I is left w.i.i .

**(d)** N(R/A) is left w.i.i .

**Theorem (2.8) :-** Let R be a Noetherian right krull homogenous ring with nilradical N and let $\Lambda$, $\Lambda'$ and A, B be as defined as in the above definitions (1.4) and (1.1) respectively. Let $P = \cap p_i$; $p_i \in \Lambda$ and let $Q = \cap q_i$; $q_i \in \Lambda'$. Then the following hold true ;

(i) N is right w.i.i ideal of R implies that N = P .

(ii) N is right w.i.i ideal of R if and only if R has an artinian quotient ring ( in this case C(P) = C(0) and hence P is a localizable semiprime ideal of R ) .

**Proof :-** (i) For the proof of (i) we first make the following claim ;

**Claim(1)** :- Let q , p be minimal prime ideals of R , with $|R/p|_r = |R|r$, and $|R/q|_r < |R|r$. Then there exists an ideal B with $|R/B|_r < |R|r$, such that $pB \leq qp$ .

Proof :- This is true because N is right w.i.i ideal of R .

**Proof of (i)** :- We now prove (i) . To prove (i) assume A1 A2 A3....Am....An = 0 , where Ai are all the minimal prime ideals of R .Then applying claim (1) we get that p1p2...pk D =0 , where each $p_i \in \Lambda$ and D is an ideal such that $|R/D|_r < |R|r$. Hence since R is right krull homogenous we must have that P is a nilpotent ideal of R . Thus N= P .

**(ii) Proof of (ii)** :- We now prove (ii) . If N is right w.i.i ideal of R , then by (i) above we have that N = P , where N is the nilradical of R . Thus using theorem (2.5) stated above we get that P is a right localizable ideal of R .

Alternatively one may use theorem (8) of [3] for the proof that R has an artinian quotient ring. Converse follows from theorem (2.5) (i) stated above.

**Corollary (2.9)** :- Let R be a Noetherian ring with nilradcal N and let A, B and Λ, Λ' be as in the definitions (1.1) and (1.3) above. Let P =∩ { Pi /Pi ε Λ } and let Q = ∩ { Qi /Qi ε Λ' } Then the following hold true ;

**(a)** N is right w.i.i if and only if N(R/A) = P/A and C(P) =C(A). Thus in this case for some intger m, $P^m \leq A$.

**(b)** N is left w.i.i if and only if N(R/B) = Q/B and C(Q) = C(B). Thus in this case for some intger m, $Q^m \leq B$.

**Proof** :- **(a)** First note that by proposition (2.6) N is right w.i.i implies that N(R/A) is right w.i.i. Since R/A is right krull homogenous, and since A ≤ P, so by theorem (2.8) above we get that N(R/A) = P /A, and hence $P^m \leq A$, for some integer m. Since from theorem (2.8) above we have that R/A has an artinian quotient ring thus we must have that c(P) = c(A). Hence by Small's theorem (see [4]) we get that P/A is a right localizable ideal of R/A. We now prove the converse. Conversely suppose for some intger m, $P^m \leq A$. Thus we have that N(R/A) = P/A .Since conversely we are given that C(P) = C(A), so by Small's theorem (see [4]) P/A is a right localisable semiprime ideal of R/A. Thus by theorem (2.5) above P/A is a right w.i.i of the ring R/A. Hence by proposition (2.6) above we get that N(R/A) is a right w.i.i ideal of R/A. Using again proposition (2.6) above we get that N is a right w.i.i ideal of R.

**(b)** The proof of (b) is on the same lines as that of (a) above.

We now state Our main theorem below ;

**Theorem (2.10)** :- Let R be a Noetherian ring with nilradical N and let $\Lambda$, $\Lambda'$, A, B be as defined in definitions (1.4) and (1.1) respectively. Let $P = \cap p_i$ ; $p_i \in \Lambda$ and $Q = \cap q_i$ ; $q_i \in \Lambda'$. Then the following hold true ;

(a) N is right w.i.i ideal of R if and only if P is a right localizable ideal of R.

**(b)** N is left w.i.i ideal of R if and only if Q is a left localizable ideal of R.

**Proof** :- **(a)** We first prove (a). Assume N is right w.i.i ideal of R. Next note that $A \leq P$. Since $|R/A|_r = |R|_r$, using theorem (2.6) stated above we get that P/A is a right w.i.i ideal of R. Hence from corollary (2.9) above we get that C(P/A) is a right ore set of the ring R/A. If $A \neq 0$, then from proposition (1.2) stated above, we have that $r(A) \neq 0$ and $|R/r(A)| < |R|$. Since P is an ideal of R such that R/P is a right krull homogenous semiprime ring with $|R/P|_r = |R|_r$, so $r(A) \not\subset P$. Thus we must have that $r(A) \cap C(P) \neq \Phi$. Let $d \in r(A) \cap C(P)$. Clearly, then $Ad = 0$. Using this together with the fact that C(P/A) is a right ore set in the ring R/A, we get that C(P) is a right ore set in the ring

R. Conversely assume that P is a right localizable ideal of R. Then use theorem (2.5) above to conclude that P is a right w.i.i of R. Hence from theorem (2.6) above we get that N is a right w.i.i ideal of R.

(b) :- The proof of (b) follows on the same lines as the proof of (a) above.

**Remark (3)** :- There is no example known of a noetherian ring R with nilradical N for which N is not a right and a left w.i.i ideal.

**Remark :- The present paper replaces an earlier paper of the author that was submitted to the arxiv without proofs.**

**Acknowledgement :- The author is much indebted to Prof. K.R.Goodearl for suggesting once briefly to him that most of the conditions on a noetherian ring that the author was considering should be viewed in terms of localizability of certain ideals of the ring.**

**References:-**

(1) K.A. Brown, T.H. Lenagan and J.T. Stafford, "Weak Ideal Invariance and localisation"; J. Lond. math.society (2), 21, 1980, 53-61.

(2) K.R. Goodaarl and R.B. Warfield ,"An Introduction to Non commutative Noetherain Rings" , L.M.S., student texts ,16.

(3) G.Krause, T.H. Lenagan and J.T. Stafford, " Ideal Invariance and Artinian Quotient Rings "; Journal of Algebra 55, 145-154 (1978).